# BandPC: Learning Residual-Band Preconditioner Combinations for Flexible Conjugate Gradient Solvers

D. M. LI [a, b, c, d] *, Shang-Tian YANG [a], Xin QIU [a, d]

a. School of Civil Engineering and Architecture, Wuhan University of Technology, Wuhan, 430070, China

b. Embodied Physical Intelligence Center in Civil Engineering, Wuhan University of Technology, Wuhan, 430070, China

c. Key Laboratory of Low-Altitude Technology and Smart Urban Renewal, Department of Housing and Urban-Rural Development of Hubei Province, Wuhan University of Technology, Wuhan, 430070, China

d. Sanya Science and Education Innovation Park of Wuhan University of Technology, Sanya, 572000, China

**Abstract:** The conjugate gradient (CG) method is a classical iterative solver for sparse symmetric positive definite (SPD) linear systems, but its convergence strongly depends on the spectral properties of the system matrix. Preconditioning can improve these properties; however, designing preconditioners that generalize across diverse and irregular sparsity patterns remains challenging. We propose BandPC (Band Preconditioner Combinations, where "Band" denotes residual-norm bands), a data-driven framework that leverages graph neural networks (GNNs) to predict multi-stage preconditioning strategies for the flexible conjugate gradient (FCG) method. BandPC partitions the FCG iteration into three residual-norm-based bands and defines a structured search space of 125 candidate combinations over five classical preconditioners. By exploiting FCG's ability to switch preconditioners across iterations, BandPC learns to map matrix structure directly to a promising preconditioner sequence. To improve training and generalization, we introduce a soft-labeling mechanism that retains near-optimal combinations and normalizes their scores into a label distribution, together with a hierarchical feature representation that captures node-level attributes, edge-level algebraic coupling strengths, and global matrix statistics. Experiments on a hybrid benchmark of synthetic SPD problems and diverse matrices from the SuiteSparse Collection show that BandPC achieves the best iteration count and solution time on 41.3% and 33.1% of test matrices, respectively, even when compared with each matrix's individually best traditional preconditioner. These results demonstrate that learned residual-band preconditioner scheduling can effectively accelerate flexible CG solvers.

*Corresponding to D. M. LI (dongmli2-c@my.cityu.edu.hk, domili@whut.edu.cn).

# 1 Introduction

Solving large-scale sparse linear systems is a fundamental task in scientific computing and industrial simulation, with applications ranging from finite element structural analysis [1] and computational fluid dynamics [2] to geophysical electromagnetic inversion [3]. Existing solution techniques fall into two broad categories: direct methods and iterative methods. Direct methods, which rely on matrix factorizations such as LU or Cholesky decomposition, offer high numerical accuracy, but their prohibitive computational complexity and memory requirements severely limit scalability in large-scale settings. Iterative methods, by contrast, refine an initial guess through successive approximations and are therefore more widely adopted in practice, especially when approximate solutions satisfy the required accuracy [4, 5].

The Conjugate Gradient (CG) algorithm is a well-established iterative solver for sparse symmetric positive definite (SPD) linear systems of the form $\boldsymbol{Ax}=\boldsymbol{b}$. However, CG exhibits slow convergence for ill-conditioned matrices, making preconditioning essential for practical performance [6].

Despite the importance of preconditioners, designing robust and generalizable ones remains a longstanding challenge [7]. Effective preconditioning must be tailored to the sparsity structure of the matrix, whereas neural network-based preconditioners are often restricted to narrow problem domains with relatively uniform sparsity patterns. This limitation motivates the development of adaptive preconditioners that can handle sparse matrices with diverse and irregular sparsity structures.

A key observation from our investigation is that different stages of the iterative convergence process impose different requirements on the preconditioner. In the early stages, when residuals are large, stronger and more aggressive preconditioners with effective smoothing properties are preferable. As convergence proceeds and residuals decrease, lightweight preconditioners become sufficient and more computationally efficient. The Flexible Conjugate Gradient (FCG) method naturally supports this idea, as it allows a different preconditioner to be used at each iteration without sacrificing convergence guarantees.

Motivated by this observation, we propose a residual-band preconditioning strategy that partitions the solving process into multiple residual-defined stages, each assigned a distinct preconditioner. This strategy systematically exploits the complementary strengths of different preconditioners along the convergence trajectory. With five candidate preconditioners distributed over three stages, it defines a structured search space of 125 candidate combinations. Because the best combination varies across matrices with different sparsity structures and exhaustive manual search is impractical, we introduce BandPC (Band Preconditioner Combinations, where “Band” denotes residual-norm bands), a Graph Neural Network (GNN)-based framework that represents a sparse matrix as a graph and efficiently predicts a promising preconditioner sequence for unseen matrices.

The main contributions of this work are summarized as follows. First, we propose a residual-band preconditioning strategy that partitions the iterative solving process into multiple residual-defined stages, each assigned a distinct preconditioner. By leveraging the dynamic switching capability of FCG, this strategy explores the synergistic effects of combining diverse preconditioners across convergence stages, yielding a structured search space of 125 candidate combinations.

Second, we introduce a soft-labeling mechanism for training data generation. Multiple near-optimal preconditioner combinations are retained, and their scores are re-normalized into a label distribution. This mitigates the overfitting issues associated with hard single-label classification in a large candidate space and improves the generalization of the trained model.

Third, we introduce a hierarchical feature representation that integrates node-level attributes, edge-level algebraic coupling strengths, and global-level statistics, enabling the GNN to capture both localized topology and macro-scale numerical characteristics.

Fourth, we propose a GNN architecture that achieves high-precision performance ranking for diverse preconditioner combinations, facilitating effective solver selection.

# 2 Related work

## 2.1 CG and Preconditioning

The Conjugate Gradient (CG) method is a well-established iterative algorithm for solving SPD linear systems. A key advantage of CG is that it avoids explicit matrix factorizations, making it efficient for large-scale sparse matrices [8]. The algorithm constructs a sequence of search directions $\{P_i\}$ from the corresponding Krylov subspace, which are mutually conjugate with respect to $A$. These directions are used to iteratively update the solution $x_k$. Since the optimal step size can be computed in closed form for each direction, the method theoretically converges within $n$ steps in exact arithmetic [9]. Typically, the initial search direction $p_0$ is chosen as the negative gradient of the quadratic functional [10].

However, a naive CG implementation often suffers from slow convergence when dealing with ill-conditioned matrices. To mitigate this, CG is typically equipped with a preconditioner $M$, which transforms the original system into $\boldsymbol{M}^{-1}\boldsymbol{Ax}=\boldsymbol{M}^{-1}\boldsymbol{b}$ to improve its spectral properties and condition number. While numerous preconditioners exist—such as Block Jacobi, Gauss-Seidel, Sparse Approximate Inverse (SPAI), and Algebraic Multigrid (AMG) [4, 11]—selecting the most effective one remains a significant challenge. In practice, practitioners often rely on a combination of theoretical heuristics and empirical trial-and-error, as the best choice is highly problem-dependent.

## 2.2 Automated Solver and Preconditioner Selection

To address this selection dilemma, recent literature has explored automated methodologies. Early approaches, such as the poly-iterative method proposed by Barrett et al. [12], executed multiple Krylov solvers simultaneously to predict the best performer. Subsequent research introduced machine learning to this domain: Bhowmick et al. [13, 14] utilized alternating decision trees for adaptive solver selection, while Holloway and Chen [15] and Kuefler and Chen [16] explored Neural Networks and Reinforcement Learning, respectively. Other ML techniques, including k-nearest neighbor (k-NN) and Support Vector Machines (SVM), have also been investigated [17]. Furthermore, specialized dynamic models have been developed for transient simulations [18], and broader algorithm selection frameworks have been evaluated by Kotthoff et al. [19]. Building on these foundational efforts, Motter et al. [20] introduced Lighthouse, an automated solver selection tool designed to simplify the complexity of selecting optimal high-performance computing (HPC) software.

As deep learning gained traction, researchers began leveraging more sophisticated neural architectures for performance prediction. Sun et al. [21] proposed a matrix feature selection strategy based on machine learning models for Krylov solver prediction, while Sappl et al. [22] demonstrated the effectiveness of Convolutional Neural Networks (CNNs) in predicting optimal preconditioners for CG solvers within the context of urban water management. More recently, to better capture the non-grid-like structure of sparse matrices, Tang et al. [23] utilized GNNs to simultaneously select both preconditioners and Krylov solvers, achieving state-of-the-art performance by learning directly from the underlying graph representation of the linear system. Xiong et al. [24] proposed MM-AutoSolver, a multimodal framework that jointly exploits numerical and structural matrix features for automated solver-preconditioner selection, and Zhou et al. [25] introduced multigrid-inspired GNNs that leverage hierarchical coarsening structures to improve selection accuracy across diverse sparse systems.

However, most of these methods treat preconditioner selection as a static, one-shot decision. They predict a single solver or preconditioner for the entire iteration and therefore do not exploit the stage-dependent behavior of preconditioners.

### 2.3 Learned Preconditioner Construction

Beyond selection, another emerging trend involves using neural networks to construct preconditioning operators. Li et al. [26] proposed learning preconditioners specifically for PDE-based CG solvers, while Häusner et al. [27] developed Neural Incomplete Factorization, which replaces traditional heuristics with learned parameters to accelerate convergence. Rudikov et al. [28] introduced the FCG-NO method, combining Flexible Conjugate Gradients with Neural Operators for efficient PDE solving. Chen [29] introduced graph neural preconditioners, further demonstrating the potential of GNNs in capturing complex topological dependencies for iterative solutions. More recently, Nastorg et al. [30] embedded a GNN preconditioner within a multi-level domain decomposition framework, extending scalability to larger problems, while Li et al. [31] combined GNNs with traditional incomplete Cholesky factorization, achieving strong generalization on sparse matrices with up to five million dimensions.

These construction methods can be powerful, but they typically learn a new preconditioning operator rather than scheduling multiple classical preconditioners. They may also incur high training costs and limited interpretability.

### 2.4 Dynamic and Multi-Stage Preconditioning

FCG permits a different preconditioner at each iteration without compromising convergence guarantees. This makes it natural to consider multi-stage preconditioning, where different preconditioners are used in different residual bands. Yet existing automated selection methods rarely address this setting. Most focus on static selection or learned operator construction, leaving the problem of learning an effective residual-band combination of multiple classical preconditioners for unseen sparse matrices largely unexplored.

BandPC fills this gap. Unlike static selection methods, it predicts a sequence of preconditioners across residual bands. Unlike learned preconditioner construction, it does not learn a new operator but instead schedules classical preconditioners. This positions BandPC as a learned multi-stage preconditioner scheduling framework for FCG.

# 3 Method

## 3.1 Problem formulation and Notation

We formulate BandPC as a supervised distribution-learning problem over candidate preconditioner combinations. Let $P = \{P_1, P_2, P_3, P_4, P_5\}$ denote five classical preconditioners. The FCG iteration is partitioned into three non-overlapping residual bands based on the residual norm $r_k = \|Ax_k - b\|_2$ at iteration $k$:

$$
\begin{aligned}
&I_1 : r_k \in [r_0, 10^{-3}], \\
&I_2 : r_k \in [10^{-3}, 10^{-6}], \\
&I_3 : r_k \in [10^{-6}, 10^{-9}].
\end{aligned}
\tag{1}
$$

For each interval $I_i$, a preconditioner $P_i$ is selected from the candidate set $P = \{None, Jacobi, IC(0), BlockJacobi, AMG\}$ .This results in a total of $|P|^3 = 125$ possible preconditioner combinations, denoted as $P = (P_1, P_2, P_3)$. Consequently, for any input matrix, the target label for our GNN model is the index of the optimal combination $i^* \in \{0,1,2,\ldots,124\}$

For each matrix in the dataset, we solve the linear system $Ax = b$ using all 125 combinations, assuming a ground-truth solution $x = 1$ for consistency. The performance of each combination $P_i$ is quantified by a score $s_i$, calculated as:

$$
s_i = \begin{cases} \log(1/t_i) & \text{if the system is solved successfully,} \\ 0 & \text{if the solver times out or fails,} \end{cases}
\tag{2}
$$

where $t_i$ represents the actual wall-clock time.

Furthermore, considering that multiple combinations may yield similar acceleration performance, we introduce a soft-labeling mechanism to enhance the model's generalization capability. We define the label distribution by retaining all combination indices that satisfy:

$$
s_i \geq 0.9 \cdot \max(\{s_j\}_{j=0}^{124})
\tag{3}
$$

The scores of these retained combinations are then re-normalized so that their sum equals 1. This refinement allows the loss function to reflect the distribution of multiple high-performance solutions, effectively mitigating the over-fitting issues typically associated with single hard labels.

### 3.2 Matrix-to-Graph Representation

We interpret a sparse matrix as the adjacency matrix of a graph. Each row/column index corresponds to a node. Off-diagonal nonzeros are mapped to edges, and diagonal entries are represented as self-loops. Since SPD matrices are symmetric, we treat the graph as undirected and merge symmetric entries. Edge weights are normalized magnitudes of the corresponding matrix entries. This representation preserves both sparsity structure and numerical coupling.

Figure 1 illustrates a representative example of this transformation, contrasting the classical matrix storage format with our proposed graph-based representation.

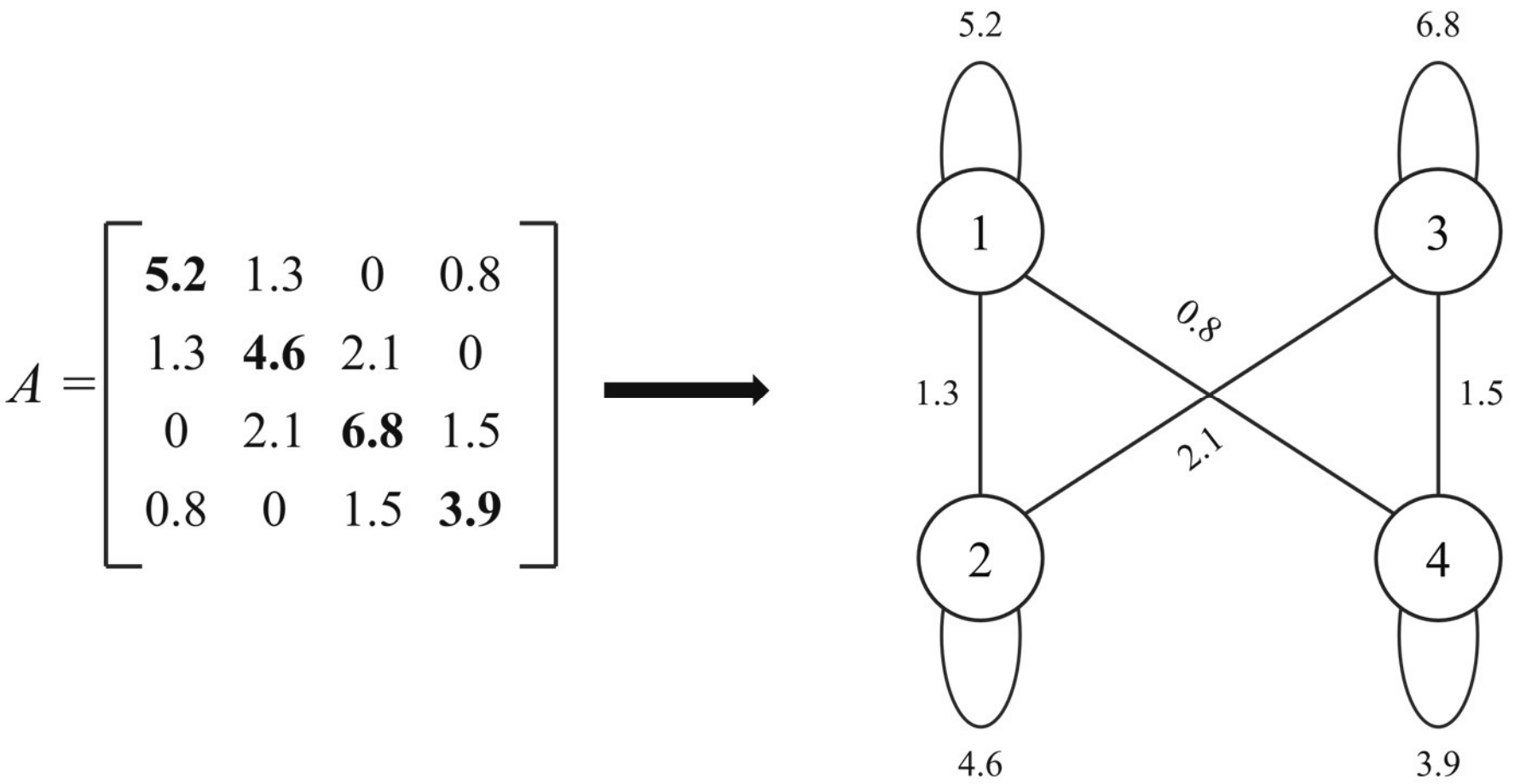

Figure 1: Illustration of the matrix-to-graph transformation.

### 3.3 Node, Edge, and Global Features

Table 1 provides a comprehensive list of the node and global features employed in this study. The local degree profile was introduced by Cai & Wang [32], while the dominance and decay features were originally proposed by Tang et al. [23]. For global features, we adopt and refine the feature set established in [20, 33] to achieve an optimal balance between classification accuracy and computational overhead. Furthermore, we optimize these features to specifically leverage the properties of SPD matrices. To ensure numerical stability when fusing heterogeneous data, all features undergo normalization or standardization prior to concatenation.

Table1： List of node and global features used in the graph neural network.

| | Features |
|---|---|
| Node | Local degree profile |
| | Diagonal dominance |
| | Diagonal decaying |
| Global | Matrix size |
| | Estimated condition number |
| | 1-norm |
| | Density |

| |
|---|
| Structure variability |
| Average distance to diagonal |
| Bandwidth |
| Number of nonzero diagonals |
| Estimated num blocks |

### 3.4 GNN Architecture

Driven by the natural and strong connection between a sparse matrix and the graph adjacency matrix, we employ a GNN to parameterize the mapping function between sparsity patterns and optimal preconditioning strategies. By treating each matrix row as a node, this architecture adaptively aggregates local connectivity information and diagonal numerical features, providing highly discriminative structural embeddings for the subsequent fusion predictor.

Figure 2 illustrates the overall workflow of the proposed GNN architecture, which processes sparse matrix inputs in Coordinate (COO) format through a dual-path feature extraction framework to simultaneously capture topological structures and numerical statistics. In the structural path, the matrix is transformed into a graph representation incorporating both explicit node features, where the latter encapsulate the numerical magnitudes of non-zero entries. This path captures complex inter-row dependencies by stacking multiple graph convolutional layers. These layers leverage a multi-head attention mechanism to adaptively weight the influence of neighboring nodes based on their structural and numerical relevance, effectively identifying algebraic coupling strengths within the matrix. Subsequently, a hybrid pooling layer combining global mean and maximum operations is utilized to compress the refined node embeddings into a fixed-dimensional structural embedding vector. Simultaneously, the statistical path extracts global features directly from the matrix and encodes them via a multilayer perceptron (MLP) to encapsulate the macro-scale context in a statistical embedding vector. These two independent paths then converge at a fusion module, where the embeddings are integrated through concatenation and fed into a deep MLP predictor. Finally, a Sigmoid activation function is applied to the output vector to generate the predicted Relative Performance Scores for the 125 preconditioner combinations.

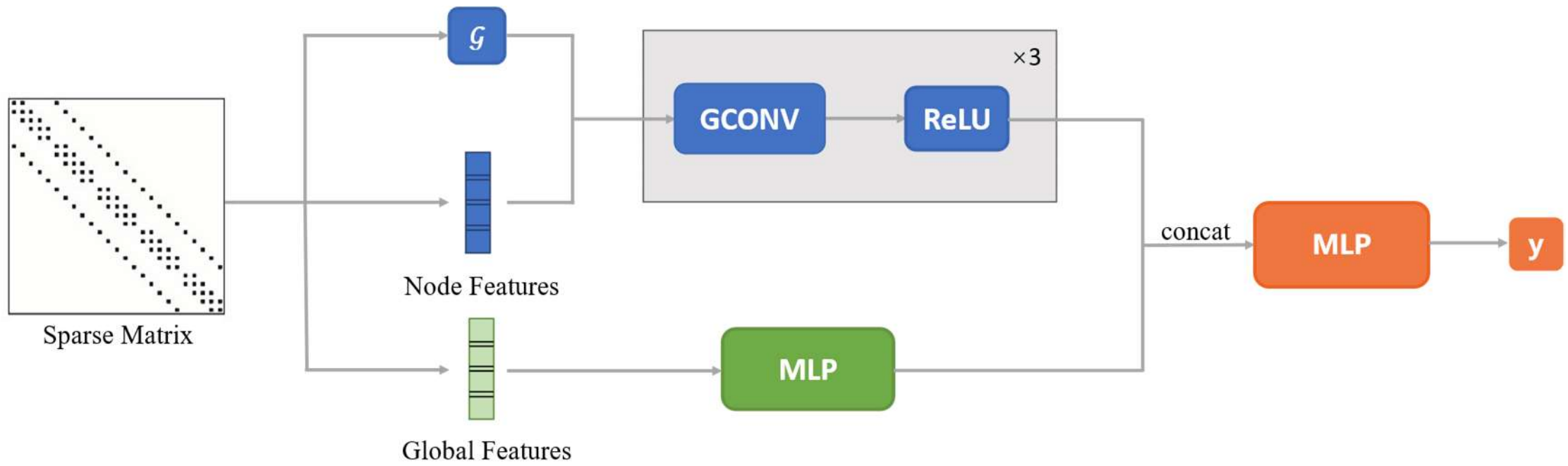

Figure 2: Overview of GNN architecture. $\mathcal{G}$ denotes the graph-level representation.

### 3.5 Soft-Label Training and Inference

The soft-label target is formed by retaining near-optimal combinations and normalizing their scores. The training loss is the binary cross-entropy (BCE) between the predicted scores and the soft-label vector:

$$\mathcal{L} = -\frac{1}{125}\sum_{i=1}^{125}\left[y_i \log(\hat{y}_i) + (1-y_i)\log(1-\hat{y}_i)\right] \tag{4}$$

where $\hat{y} \in [0,1]^{125}$. In practice, we use BCEWithLogitsLoss for numerical stability.

The model is trained with Adam, an initial learning rate of $10^{-3}$, cosine annealing, dropout of 0.3, and batch normalization. Training runs for 200 epochs with early stopping based on validation loss. At inference time, the combination with the highest predicted score is selected:

$$\hat{c}^* = \arg\max_{i\in\{1,\ldots,125\}} \hat{y}_i \tag{5}$$

The predicted combination $\hat{c}^* = (p_1^*, p_2^*, p_3^*)$ is applied across the three residual bands in FCG. Inference requires only one forward pass, and its cost is typically much smaller than the linear solve itself.

# 4 Results

## 4.1 Experimental Setup and Dataset

Our experiments were conducted on a high-performance computing server equipped with an NVIDIA GeForce RTX 3090 Ti GPU (24 GB VRAM) and a dual-processor system featuring Intel Xeon Platinum 8352V CPUs, providing a total of 144 logical threads. GNN frameworks were based mainly on PyTorch [38] and PyTorch-Geometric [39], and evaluation metrics were provided by TorchMetrics [40].

We evaluate BandPC on a hybrid dataset of 1200 sparse symmetric positive definite (SPD) matrices. This dataset consists of 116 real-world SPD matrices from the SuiteSparse Matrix Collection [41] and 1,084 synthetic SPD matrices. The SuiteSparse matrices are selected by retaining square, real-valued, SPD matrices with $1000 \le n \le 1000000$ and fewer than $2\times10^6$ nonzeros. This filtering yields 116 matrices from 21 application areas, providing diverse sparsity patterns and numerical properties.

Because the number of eligible SuiteSparse matrices is limited, training solely on them may lead to an imbalanced sample distribution and degrade generalization. We therefore generate 1,084 additional synthetic sparse SPD matrices. Each synthetic matrix $\boldsymbol{M}$ is generated from a random sparse matrix $\boldsymbol{A}$ as follows:

$$\boldsymbol{M} = \boldsymbol{A}\boldsymbol{A}^{\mathrm{T}} + \alpha\mathbf{I} \tag{6}$$

where $\alpha \approx 10^{-3}$ is used to ensure that $\boldsymbol{M}$ is symmetric positive definite. The resulting matrices are controlled to remain within the size and nonzero ranges used in our experiments. To avoid duplicate or overlapping instances, we use unique and non-overlapping random seeds for synthetic matrix generation.

The 1200 matrices are split into 80% training and 20% testing, corresponding to 960 training matrices and 240 test matrices. The split is performed at the matrix level and is stratified to preserve the ratio between synthetic and SuiteSparse matrices, avoiding leakage of the same problem family or generation seed across splits. To support hyperparameter selection and early stopping, we further separate 10% of the training set for validation, resulting in 864 training, 96 validation, and 240 test matrices. Final performance is reported on the independent test set.

Our evaluation strives to encompass a diverse range of problems across multidisciplinary domains. To this end, we utilize the SuiteSparse Matrix Collection [41], a prominent benchmark in numerical linear algebra. We specifically select square, real-valued, and SPD matrices whose number of rows falls between 1K and 1M and whose number of nonzeros is fewer than 2M. This selection results in 116 matrices from 21 application areas. However, the limited number of eligible matrices from the repository poses a challenge for achieving robust statistical representativeness and ensuring strong generalization. To mitigate the risk of overfitting on a constrained sample pool, we supplement the SuiteSparse matrices with a dataset of synthetic problems. By integrating these synthetic instances with SuiteSparse matrices, we construct a more comprehensive and balanced benchmark for training and evaluating our model.

Table 2 summarizes the performance of standard baseline preconditioners across both datasets, reporting metrics only for matrices that achieved convergence. The stopping criteria of all linear system solutions are *rtol* = $10^{-9}$ and *maxiters* = 100000.

Table2： Mean results for the synthetic and SuiteSparse dataset. $t_b$, $t_s$, and $t$ correspond to building time, solving time, and total time of the linear system.

| Dataset | Preconditioner | Iterations | $t_b$ | $t_s$ | $t$ |
|---|---|---|---|---|---|
| Synthetic | None | 2582.83 | - | 30.228 | 30.228 |
| | Jacobi | 1338.27 | 0.004 | 7.803 | 7.807 |
| | BlockJacobi | 1932.68 | 0.022 | 16.373 | 16.395 |
| | IC(0) | 103.08 | 0.367 | 0.516 | 0.883 |
| | AMG | 234.31 | 0.104 | 6.556 | 6.660 |
| SuiteSparse | None | 9425.38 | - | 60.957 | 60.957 |
| | Jacobi | 3825.45 | 0.005 | 9.604 | 9.609 |
| | BlockJacobi | 3050.05 | 0.019 | 19.917 | 19.936 |
| | IC(0) | 14.47 | 0.290 | 0.089 | 0.379 |
| | AMG | 895.04 | 0.123 | 15.511 | 15.634 |

## 4.2 Acceleration of Predicted Combinations

We compared the optimal combinations predicted by our proposed strategy with traditional preconditioning methods. The experimental results demonstrate that the optimal combinations identified by our model exhibit significant acceleration advantages in the test scenarios.

Figure 3 demonstrates representative convergence histories where our proposed method, BandPC, achieves superior performance across diverse problem domains.

It is observed that the convergence behavior of standard preconditioners varies significantly depending on the underlying matrix structure. For instance, IC(0) fails to converge on the *bcsstk36* matrix. Excluding BandPC, the most efficient baseline strategies for *bcsstk36*, *crankseg_1*, and *LF10000* are Block Jacobi, IC(0), and None (no preconditioning), respectively.

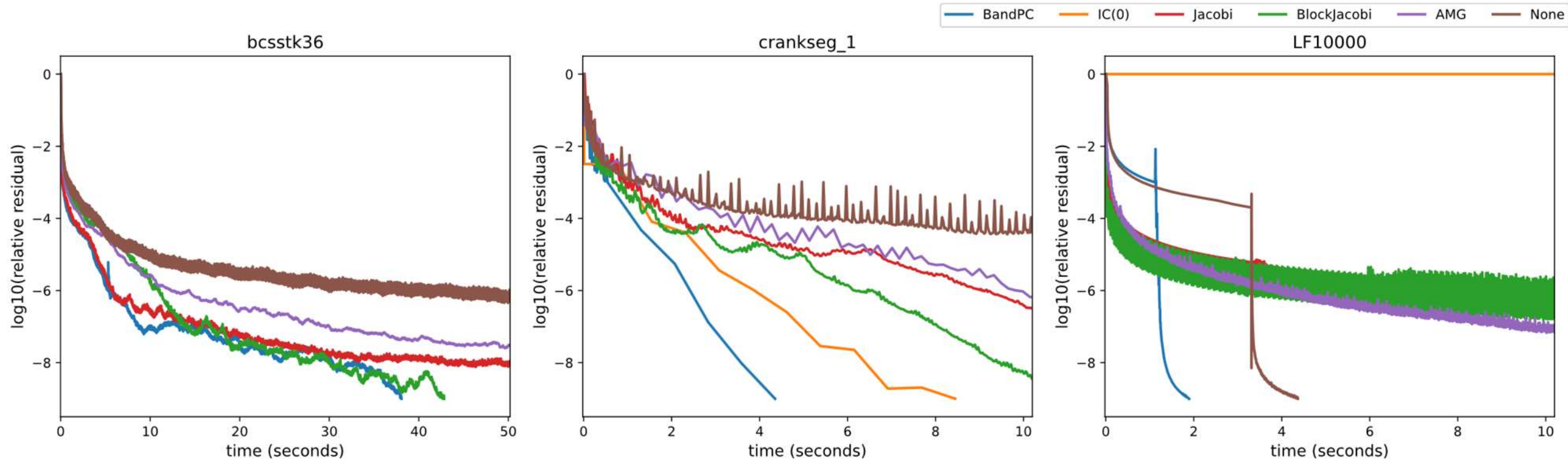

Figure 3: The convergence histories of examples.

In contrast, BandPC (indicated by the blue curve) reaches the target relative residual ($10^{-9}$) significantly faster. The specific optimal combinations identified by our model are as follows:

- For *bcsstk36*: P =（Jacobi, IC(0), IC(0)）
- For *crankseg_1*: P =（None, Block Jacobi, None）
- For *LF10000*: P =（Jacobi, Jacobi, Block Jacobi）

Compared to traditional single-level preconditioners, BandPC effectively overcomes the numerical stiffness that often leads to convergence stagnation. A striking example is seen in the *LF10000* case, where standard baselines encounter a performance bottleneck at a residual level of approximately $10^{-3}$. Conversely, the strategy prioritized by our model demonstrates exceptional numerical robustness, maintaining a steep descent toward convergence. This highlights BandPC's capability to adaptively navigate the complex search space of preconditioning sequences to find the most efficient solution for challenging linear systems.

It is worth noting that BandPC naturally simplifies to a traditional single-stage preconditioner if one specific method is identified as the optimal choice across all three stages. This inherent flexibility ensures that our model remains compatible with classical strategies while possessing the capacity for multi-stage innovation. Comprehensive convergence curves, evaluated in terms of both execution time and iteration counts.

To quantify the performance of various preconditioning strategies, we adopt two evaluation metrics from [29]: Iter-AUC and Time-AUC. These metrics represent the area under the convergence curve, providing a holistic measure of both convergence efficiency and computational speed.

The first metric, which measures the area under the relative residual norm curve with respect to iterations, is defined as:

$$\text{Iter-AUC} = \sum_{i=0}^{\text{iters}} \log_{10} r_i - \log_{10} \text{rtol}, \quad r_i = \| \mathbf{b} - \mathbf{A}\mathbf{x}_i \|_2 / \| \mathbf{b} \|_2, \tag{7}$$

where iters denotes the total number of iterations required for the solver to satisfy the stopping criterion and $r_i$ is the relative residual norm at iteration $i$. A smaller Iter-AUC value indicates superior error-reduction efficiency per iteration.

The second metric evaluates the area under the relative residual norm curve with respect to elapsed time, defined as:

$$\text{Time-AUC} = \int_0^T [\log_{10} r(t) - \log_{10} \text{rtol}] dt \approx \sum_{i=1}^{\text{iters}} [\log_{10} r_i - \log_{10} \text{rtol}](t_i - t_{i-1}), \tag{8}$$

where $T$ is the total elapsed time when the solver stops, $r(t)$ is the relative residual norm at time $t$, and $t_i$ is the elapsed time at iteration $i$. While Iter-AUC focuses on mathematical convergence, Time-AUC accounts for the actual computational overhead, making it a critical indicator for practical solver performance.

In Figure 4 we show the percentage of problems on which each preconditioner performs the best, with respect to iteration counts and solution time.

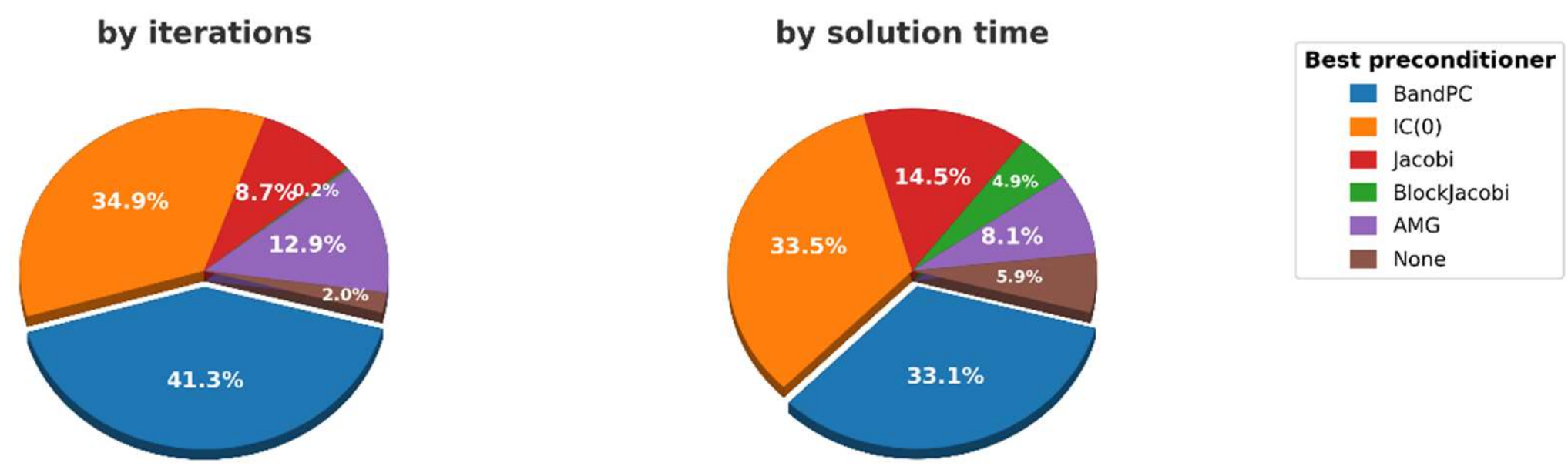

Figure 4: Percentage of problems on which each preconditioner performs the best.

It is worth noting that in the pie charts presented, BandPC refers exclusively to mixed preconditioning combinations that differ from any single traditional preconditioner applied uniformly across all intervals. Therefore, these figures demonstrate that even under the assumption that every matrix is assigned its individually optimal traditional preconditioner, BandPC still achieves the best performance in terms of both iterations (41.3%) and solution time (33.1%) on a substantial portion of the test matrices.

### 4.3 GNN Prediction Quality

While Section 4.2 demonstrated the significant speedup achieved by the optimal BandPC configurations, such combinations cannot be determined a priori without actual computation. Therefore, the efficacy of BandPC is fundamentally contingent upon the GNN model's ability to learn representative matrix features and prioritize high-performing solver sequences.

Our neural network is built upon Graph Attention Networks (GATs), which are well-suited for our task due to their native support for edge features. Specifically, edge features are projected via learnable transformation matrices and subsequently concatenated with the corresponding node features, enabling the model to jointly leverage both structural and value-level information of the sparse matrix.

We adopt two standard evaluation metrics for multilabel classification: Label Ranking Average Precision(LRAP) [30] and Normalized Discounted Cumulative Gain (NDCG) [31-33]. Both metrics operate directly on predicted probability scores without requiring a decision threshold, and return values in [0, 1], where higher scores indicate better ranking quality. Furthermore, we introduce a novel evaluation metric, the Speedup Retention Rate (SRR), to more directly characterize the practical acceleration efficacy of the proposed method.

SRR is defined as the ratio of the solving time achieved by the GNN-predicted combination to that of the ground-truth optimal combination, formally measuring how much of the theoretically achievable speedup is retained by the model's predictions. A SRR value close to 1 indicates that the predicted combination achieves performance nearly identical to the optimal, demonstrating the practical effectiveness of the proposed approach beyond what standard classification metrics can capture.

LRAP measures the average precision of the predicted ranking by evaluating the fraction of relevant labels ranked above each relevant label. A higher LRAP score indicates better alignment between the predicted scores and the true label distribution. The definition of LRAP is as follows:

$$LRAP(y,\hat{y})=\frac{1}{n}\sum_{i=0}^{n-1}\frac{1}{\|y_i\|_0}\sum_{j:y_{ij}=1}\frac{\{k:y_{ik}=1,\hat{y}_{ik}>\hat{y}_{ij}\}}{rank_{ij}} \tag{9}$$

We denote the true labels by $y$ and the prediction by $\hat{y}$ . where $n$ is the number of samples, $rank_{ij}=\left|\{k:\hat{y}_{ik}>\hat{y}_{ij}\}\right|$, and $\|\cdot\|_0$ computes the $l_0$-norm, which is equal to the number of nonzeros.

NDCG evaluates ranking quality with a position-aware discount, assigning greater importance to correctly ranked combinations at the top. It is based on the DCG score, defined as:

$$DCG(y,\hat{y})=\sum_{r=1}^{m}\frac{y_{f(r)}}{\log(1+r)} \tag{10}$$

where $m$ is the number of classes and $f$ is a ranking function induced from $y$ and $\hat{y}$. The NDCG score is the DCG score divided by the DCG score of $y$, which is considered as the ideal DCG score.

SRR directly quantifies the practical acceleration efficacy of the model's predictions, measuring how much of the theoretically achievable speedup is retained when the GNN-predicted combination is used instead of the ground-truth optimal. It is defined as:

$$SRR=\frac{t_{pred}}{t_{opt}} \tag{11}$$

where $t_{pred}$ is the solving time of the predicted combination and $t_{opt}$ is that of the ground-truth optimal combination.

Table 3 LRAP, NDCG and SRR scores of the considered methods on the test dataset.

Table3：Classification scores of considered methods on train and test dataset.

| Metrics | train | test |
|---|---|---|
| LRAP | 0.9266 | 0.9136 |
| NDCG | 0.8307 | 0.8012 |
| SRR | 0.9283 | 0.8908 |

In Figure 5, we compare the solution time of each matrix in the test set across different preconditioners. Here, BandPC denotes the preconditioning sequence predicted by our GNN model, rather than the theoretical optimum. Using the same stopping criteria defined in Section 4.1, we observe that BandPC consistently ranks among the top two in computational efficiency for the clear majority of cases, with no instances of significant performance degradation.

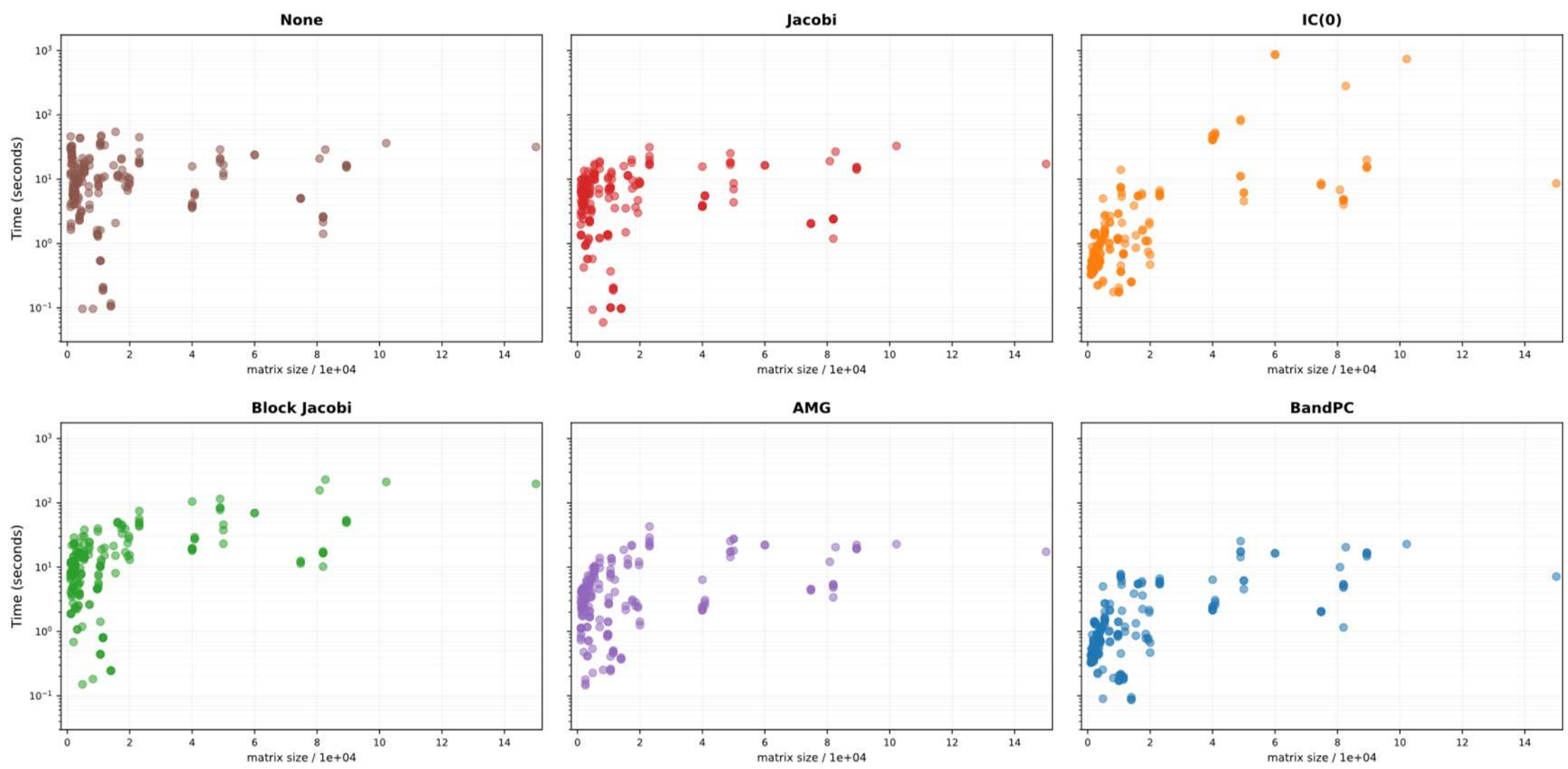

Figure 5: Preconditioner solution time.

## 5 Conclusions

In this paper, we proposed BandPC, a GNN-based framework for learning residual-band preconditioner combinations in flexible conjugate gradient solvers. BandPC partitions the FCG iteration into three residual bands and predicts a promising preconditioner sequence from 125 candidate combinations over five classical preconditioners. It further introduces soft labeling, hierarchical feature representation, and a GNN ranking architecture to improve generalization across diverse sparsity patterns.

Experiments on synthetic SPD problems and SuiteSparse matrices demonstrate that BandPC achieves the best iteration count and solution time on 41.3% and 33.1% of test matrices, respectively, even when compared with each matrix's individually best traditional preconditioner. These results indicate that learned multi-stage preconditioner scheduling can effectively accelerate flexible CG solvers.

A limitation of the current design is that it targets SPD systems, uses a fixed residual-band partition, and relies on exhaustive label generation over a relatively small candidate set. Future work will extend BandPC to non-SPD systems, incorporate preconditioner construction and switching costs into the objective, explore adaptive band boundaries, and investigate more efficient label generation and larger preconditioner libraries.

**Declaration of Generative AI and AI-assisted technologies**

During the preparation of this work, the authors used Deepseek to check grammar and improve sentence clarity. The tool was not used to generate data, analyse results, or produce scientific conclusions. All content was reviewed and verified by the authors, who take full responsibility for the accuracy and integrity of the published work.